\documentclass[oneside]{amsart}
\def\homepage#1{\urladdr{\url{#1}}}
\def\volumeyear#1{}
\def\classification#1{\subjclass[2000]{#1}}

\title{Azumaya Objects in Triangulated Bicategories}
\date{March 19, 2013}

\author{Niles Johnson}

\homepage{http://www.nilesjohnson.net}

\email{niles@math.ohio-state.edu}

\address{Dept. of Mathematics\\
The Ohio State University Newark\\
1179 University Drive\\
Newark, OH 43055 USA}

\volumeyear{2013}

\classification{55U99,18D35,16K50,14F22}

\usepackage[T1]{fontenc}
\usepackage[english]{babel}
\usepackage{stackrel} 

\usepackage{geometry}

\usepackage[unicode=true, pdfusetitle,
 bookmarks=true,bookmarksnumbered=false,
 breaklinks=false,
 backref=false,
 final
]{hyperref}

\usepackage{bookmark}

\usepackage{ifpdf}
\ifpdf
\else
\PassOptionsToPackage{dvistyle}{todonotes}
\fi

\usepackage[capitalise]{cleveref}
\newcommand{\fref}{\cref}
\newcommand{\Fref}{\Cref}

\crefname{lem}{Lemma}{Lemmas}
\crefname{thm}{Theorem}{Theorems}
\crefname{defn}{Definition}{Definitions}
\crefname{prop}{Proposition}{Propositions}
\crefname{rmk}{Remark}{Remarks}
\crefname{cor}{Corollary}{Corollaries}
\crefname{figure}{Figure}{Figures}

\usepackage{amsmath,amsthm}

\numberwithin{equation}{section} 
\numberwithin{figure}{section} 

\theoremstyle{plain} 
\newtheorem{thm}{Theorem}[section]
\newtheorem*{thm*}{Theorem}
\newtheorem{prop}[thm]{Proposition}
\newtheorem{lem}[thm]{Lemma}
\newtheorem{cor}[thm]{Corollary}

\theoremstyle{definition} 
\newtheorem{defn}[thm]{Definition}
\newtheorem{defns}[thm]{Definitions}
\newtheorem{example}[thm]{Example}
\newtheorem{notn}[thm]{Notation}

\theoremstyle{remark} 
\newtheorem{rmk}[thm]{Remark}

\newtheorem*{note}{Note}

\usepackage[all,cmtip]{xy}
\xyoption{curve}      

\usepackage{amsmath,amssymb}

\usepackage{eucal}

\DeclareMathAlphabet{\mathscr}{OT1}{pzc}%
                                 {m}{it}

\DeclareFontFamily{OMS}{rsfs}{\skewchar\font'60}
\DeclareFontShape{OMS}{rsfs}{m}{n}{<-5>rsfs5 <5-7>rsfs7 <7->rsfs10 }{}
\DeclareSymbolFont{rsfs}{OMS}{rsfs}{m}{n}
\DeclareSymbolFontAlphabet{\scr}{rsfs}

\newcommand{\sB}{\scr{B}}
\newcommand{\sC}{\scr{C}}
\newcommand{\sD}{\scr{D}}
\newcommand{\sE}{\scr{E}}

\newcommand{\sM}{\scr{M}}

\newcommand{\sO}{\scr{O}}
\newcommand{\sP}{\scr{P}}

\newcommand{\sR}{\scr{R}}
\newcommand{\sS}{\scr{S}}
\newcommand{\sT}{\scr{T}}

\newcommand{\bP}{\mathbb{P}}

\newcommand{\bZ}{\mathbb{Z}}

\newcommand{\al}{\alpha}
\newcommand{\be}{\beta}

\newcommand{\epz}{\varepsilon}

\newcommand{\Ga}{\Gamma}

\newcommand{\Si}{\Sigma}

\newcommand{\id}{\text{id}}

\newcommand{\fto}{\xrightarrow}

\newcommand{\cn}{\colon}

\newcommand{\fl}[1]{\lfloor #1 \rfloor}
\newcommand{\lng}{\left\langle}
\newcommand{\rng}{\right\rangle}
\newcommand{\wt}{\widetilde}

\newcommand{\cmxymat}{\xymatrix}

\newcommand{\cell}[3]{#1:#2 \sto #3}

\usepackage{slashed}
\declareslashed{}{\text{\rule[1.5pt]{.2pt}{3pt}}}{-.05}{0}{\rightarrow}
\newcommand{\sto}{\slashed{\rightarrow}}
\newcommand{\sfto}[1]{{%
  \declareslashed{}{\text{\rule[-9pt]{.2pt}{3pt}}}{-0.05}{0}{\xrightarrow{#1}}\slashed{\xrightarrow{#1}}}}

\DeclareMathOperator{\Mon}{Mon}
\DeclareMathOperator{\End}{End}

\DeclareMathOperator{\sHom}{\textsl{s}Hom}
\DeclareMathOperator{\tHom}{\textsl{t}Hom}

\newcommand{\shom}[2]{\sHom(#1,#2)}
\newcommand{\thom}[2]{\tHom(#1,#2)}

\newcommand{\thomshort}[5]{\thom{#2}{#4}}

\newcommand{\yda}[2][\sY]{#1(#2,-)} 
\newcommand{\oyda}[2][\sY]{#1(-,#2)}

\usepackage{mathrsfs}
\usepackage{stmaryrd}

\DeclareMathAlphabet{\mathpzc}{OT1}{pzc}{m}{it}
\DeclareMathOperator{\Hom}{Hom}

\newcommand{\Ex}{\scr{E}x}
\newcommand{\Roid}{\sR\textit{ingoids}}
\newcommand{\V}{\textbf{V}}
\newcommand{\Ab}{\mathpzc{Ab}}

\newcommand{\op}{\mathsf{op}}

\newcommand{\hty}{\simeq}
\newcommand{\iso}{\cong}

\newcommand{\sma}{\wedge}

\newcommand{\aA}{{\CMcal A}}

\newcommand{\pop}{\,{\square}\,}

\begin{document}

\begin{abstract}
  We introduce the notion of Azumaya object in general
  homotopy-theoretic settings.  We give a self-contained account of
  Azumaya objects and Brauer groups in bicategorical contexts,
  generalizing the Brauer group of a commutative ring.  We go on to
  describe triangulated bicategories and prove a characterization
  theorem for Azumaya objects therein.  This theory applies to give a
  homotopical Brauer group for derived categories of rings and ring
  spectra.  We show that the homotopical Brauer group of an
  Eilenberg-Mac~Lane spectrum is isomorphic to the homotopical Brauer
  group of its underlying commutative ring.  We also discuss tilting
  theory as an application of invertibility in triangulated
  bicategories.
\end{abstract}

\maketitle

\section{Introduction}\label{sec:intro}

The notion of Azumaya algebra over a commutative ring $k$ generalizes
the notion of a central simple algebra over a field.  This was
introduced by Azumaya \cite{Azu51maximally} in the case that $k$ is
local and by Auslander-Goldman \cite{AG60Brauer} for general $k$.  The
main classification theorem of Azumaya algebras says that the set of
Morita-equivalence-classes of Azumaya $k$-algebras is the maximal
subgroup of the monoid formed by Morita-equivalence-classes of
$k$-algebras under $\otimes_k$.  This is the Brauer group of $k$, an
invariant which carries interesting algebraic and geometric
information.

\subsection{Foundations of homotopical Brauer theory}

This work aims to develop foundations for Brauer theory in homotopical
settings.  We consider Azumaya objects in closed autonomous symmetric
monoidal bicategories, and in particular focus on the triangulated
bicategories arising as homotopy bicategories of rings and ring
spectra.  (The term ``closed'' refers to the existence of internal
homs, and ``autonomous'' refers to opposites such as opposite
algebras. See \fref{sec:Closed-Str,sec:auto-str}.)

For the presentation here, we have three audiences in mind: the
algebraic audience, for whom Azumaya algebras and Brauer groups of
commutative rings are quite familiar; the categorical audience, for
whom bicategories and invertibility therein are quite familiar; and
the topological audience for whom structured ring spectra and the
homotopical algebra thereof are quite familiar.  We expect few readers
to be in the intersection of these three audiences, and thus have
attempted to write for the union of their complements.

In \fref{sec:intro,sec:Mot-Exa} we introduce and motivate our results
in the discrete and homotopical cases with a minimum of bicategorical
language.  This includes a review of classical
Brauer groups for commutative rings as a special case of our general
approach.  A survey of the relevant bicategorical background is then
presented in \fref{sec:context}, with the general characterization of
Azumaya objects in \fref{sec:AzumayaBrauer}.  We develop the notion of
Azumaya object for triangulated bicategories in \fref{sec:trbicat},
and discuss homotopical (derived) applications in
\fref{sec:homotopical-brauer-groups}.
\\

Now we state the main definitions and results.

\begin{defn}[Eilenberg-Watts Equivalence]
  Let $A$ and $B$ be a 0-cells of a bicategory $\sB$.  We say $A$ is
  \emph{Eilenberg-Watts equivalent} to $B$ if there exists an
  invertible 1-cell $\cell{T}{A}{B}$.  In bicategorical literature,
  this is known simply as \emph{equivalence} of 0-cells, but we prefer
  the more expressive term for this work.
\end{defn}

\begin{defn}[Eilenberg-Watts Bicategory]
  We say that a bicategory is \emph{Eilenberg-Watts} if it is closed
  (i.e., has internal homs); has a symmetric monoidal product; and has
  an autonomous structure (i.e., a notion of ``opposite'' objects).
  Each of these bicategorical structures is described in
  \cref{sec:context}.
\end{defn}

\begin{defn}[Brauer Group, Azumaya Objects]
  Let $\sB$ be a monoidal bicategory with unit 0-cell $k$.  The
  \emph{Brauer group} of $\sB$, denoted $Br(\sB)$, is the group of
  1-cell equivalence classes (Eilenberg-Watts equivalence classes) of
  0-cells $A$ for which there exists a 0-cell $B$ such that $A
  \otimes_k B$ is Eilenberg-Watts equivalent to $k$.  Such 0-cells are
  called \emph{Azumaya objects}.
\end{defn}

\begin{rmk}
  What we have termed ``Eilenberg-Watts equivalence'' is sometimes
  called ``Morita equivalence'', or sometimes ``standard Morita
  equivalence'' \cite{konig1998deg}.  The term ``Morita equivalence''
  is also used frequently for monoid objects to mean that their
  categories of modules are equivalent.  The fundamental theorem of
  Morita theory states that these notions coincide for algebras over a
  commutative ring and their module categories.  Moreover, Rickard
  \cite{rickard1989mtd} has shown that these notions coincide for
  algebras over a commutative ring and their derived categories of
  modules.  However it is known that this coincidence does not
  generalize to derived categories of differential graded algebras or
  ring spectra \cite{shipley2006mts}.  Work of Dugger and
  Dugger--Shipley analyzes the situation for differential graded
  algebras and ring spectra; see \cite{DS07Topological,DS07Enriched}
  and their related articles.  The survey \cite{johnson2008mtd}
  compares the two notions abstractly for closed bicategories.

  We have elected to introduce new terminology in hopes of clarifying
  the distinction between the two notions.  See
  \cite{Hov09Eilenberg-Watts} for a discussion of Eilenberg-Watts
  theorems in various settings.
\end{rmk}

The classical treatment of Azumaya algebras depends heavily on
localization arguments and reducing to the case of algebras over
fields; one feature of recent work in this area is a treatment that is
independent of ideal theory:  by introducing a bicategorical context
we are able to generalize and unify the basic theory of Azumaya
algebras and Brauer groups, giving elementary formal proofs of the
main classification theorems.  Similar ideas have appeared in
algebraic and categorical literature---notably
\cite{Par75Non-Additive,Vit96Brauer,VOZ98Brauer,borceux2002ac}.  We
extend this approach to homotopical settings, giving a classification
of Azumaya objects via invertibility in triangulated bicategories in
\fref{sec:invertibility-trbicat}.  During preparation of this article,
related work has appeared in \cite{BRS10Brauer,To12Derived,AG12Brauer,DAT12Quillen}.

\begin{notn}
  Throughout, $k$ will denote the unit of a symmetric monoidal
  bicategory.  In our applications, $k$ will often be a commutative
  ($E_\infty$) d.g.~algebra or ring spectrum.  We let $\sC_k$ denote
  the category of $k$-modules, noting that this is a symmetric
  monoidal category.  We let $\sM_k$ denote the bicategory of
  $k$-algebras and bimodules.

  When $\sC_k$ has a suitable model structure (described in
  \fref{sec:hty-bicat}, following \cite{SS00Algebras}), then there is
  an induced model structure on $\sM_k(A,B)$ for $k$-algebras $A$ and
  $B$.  We let $\sD_k(A,B)$ denote the corresponding homotopy
  category.  

  Taking as 0-cells the collection of $k$-algebras which are cofibrant
  as $k$-modules, we give one set of conditions in
  \fref{prop:model-bimod-bicat} which guarantee that $\sD_k$ forms a
  triangulated bicategory with the additional Eilenberg-Watts
  structure necessary for the homotopical Brauer theory of
  \fref{sec:AzumayaBrauer,sec:invertibility-trbicat}.
\end{notn}

Although the bicategorical definitions of Azumaya object and Brauer
group require only a monoidal structure, these are somewhat too
general to be of practical interest.  In an Eilenberg-Watts
bicategory, we have characterization theorems for Azumaya objects in
more concrete terms.  In
\cref{defn:general-central-endomorphic,defn:general-fproj-fsep} we
generalize the classical notions of central, (faithfully) separable,
and faithfully projective $k$-algebras to Eilenberg-Watts
bicategories.  We also introduce the term ``endomorphic'', formally
dual to the notion of centrality.  We then show that these are the
correct generalizations by proving the following theorem in
\fref{sec:AzumayaBrauer}:

\begin{thm}\label{thm:AzumayaEquiv} Let $A$ be a 0-cell of an
  Eilenberg-Watts bicategory, and let
    \[
    \cell{\fl{A}}{k}{A^e} = A \otimes A^{\op}
    \]
    denote the 1-cell induced by $A$ and the autonomous structure.
    The following are equivalent:
  \begin{enumerate}
  \item\label{AzumayaEquivInv} $\fl A$ is an invertible 1-cell.
  \item\label{AzumayaEquivRtEC} $A$ is central and faithfully separable over
    $k$.
  \item\label{AzumayaEquivLtEC} $A$ is endomorphic and faithfully projective
    over $k$.
  \item $A^e$ is Eilenberg-Watts equivalent to $k$.
  \item There is a 0-cell $B$ such that $A \otimes B$ is Eilenberg-Watts equivalent to
    $k$.
  \end{enumerate}
\end{thm}
\begin{rmk}
  The equivalence of the first three conditions is formal and
  straightforward.  The argument proceeds by explaining that each of the
  second two conditions is an alternate description of invertibility.  It
  is clear that the first condition implies the last two; one goal of this
  work is to develop a setting in which the reverse implications are also
  straightforward.
\end{rmk}


Another goal of this work is to prove a further characterization of
Azumaya objects when (as in the derived and homotopical cases) the
ambient bicategory carries a triangulated structure.  In this case,
the notion of invertibility admits alternate descriptions in
terms of localization.  This has been used in the topological case by
Baker-Lazarev \cite{baker2004thc} for $THH$ calculations.  In
\cref{sec:trbicat} we introduce triangulated bicategories and prove a
further characterization of Azumaya objects in that context.  As with
the rest of this theory, the notion of localization has both left
(source) and right (target) versions, but we leave these definitions
until \cref{sec:trbicat} where we will develop more of the underlying
bicategorical structure.  The following is immediate from
\fref{thm:InvertibilityChar}:

\begin{thm}
  \label{thm:Azumaya-char-triang}
  Let $A$ be a 0-cell in a triangulated Eilenberg-Watts bicategory.
  The following statements are equivalent:
\begin{enumerate}
  \item $A$ is an Azumaya object.
  \item $A$ is dualizable over $A^e$ and central over $k$, and $A^e$ is
    strongly (target-)$A$-local.
 \item $A$ is dualizable over $k$ and endomorphic over $k$, and $k$ is
    strongly (source-)$A$-local.
 \end{enumerate}
\end{thm}

Using these general characterizations of Azumaya objects, we show in
\fref{sec:homotopical-brauer-groups} that the homotopical Brauer
groups defined here agree with those of \cite{To12Derived} for derived
categories of rings and of \cite{BRS10Brauer} for commutative
$S$-algebras.


\subsection{Brauer groups of Eilenberg-Mac~Lane spectra}

One application of the foundations developed here is an immediate
computation of the Brauer group of an Eilenberg-Mac~Lane spectrum.  We
include the proof because it is short and because it illustrates the
utility of the setting we develop: that which should be formal
(because the homotopical algebra of Eilenberg-Mac~Lane spectra is
equivalent to the ordinary algebra of their underlying rings) has
become straightforward (because the Eilenberg-Mac~Lane functor induces
a local equivalence of triangulated bicategories).

\begin{thm}\label{thm:em-brauer-iso}
  Let $k$ be a commutative ring.  The Eilenberg-Mac~Lane functor to
  symmetric spectra
  \[
  H\cn \sC_k \to \sC_{Hk}
  \]
  induces an isomorphism of Brauer groups
  \[
  H\cn Br(\sD_k) \to Br(\sD_{Hk}).
  \]
\end{thm}
\begin{proof}
  By \fref{rmk:k-alg-trbicat,prop:homotopy-bicat-applications} we have triangulated
  Eilenberg-Watts bicategories $\sD_k$ and $\sD_{Hk}$.  The
  Eilenberg-Mac~Lane functor $H$ is a weak monoidal Quillen
  equivalence, which means that $H$ is a lax monoidal functor (but its
  adjoint is not necessarily monoidal).  Thus $H$ sends differential
  graded $k$-algebras to $Hk$-algebras, and sends $(A,B)$-bimodules to
  $(HA,HB)$-bimodules \cite{Shi07HZ-algebra}.  We have natural
  equivalences in $\sD_{Hk}(HC,HB)$
  \[
  H(X \sma_A W) \hty HX \sma_{HA} HW
  \]
  for $X \in \sD_k(C,A)$ and $Y \in \sD_k(A,B)$.  Therefore $H$
  induces a pseudofunctor
  \[
  \sD_k \to \sD_{Hk}.
  \]
  Moreover, it gives an equivalence
  \[
  \sD_k(A,B) \fto{\hty} \sD_{Hk}(HA, HB)
  \]
  for all differential graded $k$-algebras $A$ and $B$ and hence a
  bijection on isomorphism classes of invertible bimodules and on
  Eilenberg-Watts equivalence classes of Azumaya objects.  Thus we
  have an isomorphism of Brauer groups.
\end{proof}

\begin{rmk}
  Note that this result depends on the base commutative ring spectrum
  $Hk$ being Eilenberg-Mac~Lane.  If $Hk$ is an algebra over another
  commutative ring spectrum $R$, then restriction along the unit map
  $R \to Hk$ yields a pseudofunctor of bicategories and hence a map of
  Brauer groups
  \[
  Br(\sD_k) \fto{\iso} Br(\sD_{Hk}) \to Br(\sD_R)
  \]
  which may be neither injective nor surjective.  The interesting
  example of \cite[\S 5]{shipley2006mts} shows that restriction from
  $H\bZ$ to the sphere spectrum is non-injective on
  Eilenberg-Watts equivalence classes of $H\bZ$-algebras.
\end{rmk}

As an immediate corollary, we have a computation of the Brauer group
for an Eilenberg-Mac~Lane spectrum of a field in terms of the
classical Brauer group:
\begin{cor}\label{thm:brauer-gp-field-iso}
  If $k$ is a field, then $Br(\sD_{Hk}) \iso Br(k)$.  In particular, $Br(\sD_{Hk})
  = 0$ if $k$ is finite or algebraically closed.
\end{cor}
\begin{proof}
  \fref{prop:derived-brauer-group-field} says that $Br(\sD_k) \iso Br(k)$.
  The result then follows from \fref{thm:em-brauer-iso}.
\end{proof}
\noindent In the case that $k$ is an algebraically closed field, this
was established by \cite{BRS10Brauer} using different methods, but
also relying on \cite{To12Derived} for the comparison with the derived
Brauer group of a field.

\subsection{Invertibility}

\cref{sec:dltybicat,sec:invertibility-trbicat} give a general study of
invertibility in triangulated bicategories, and this forms the
foundation of the characterization theorems.
\cref{thm:InvertibilityChar} gives a general characterization of
invertible objects in triangulated bicategories.  As further
applications, we recover results from the derived Morita
(Eilenberg-Watts) theory of Rickard \cite{rickard1989mtd} and
Schwede-Shipley \cite{SS03Stable}.  The object $T$ here is called a
\emph{tilting object}.

\begin{prop}[\cite{rickard1989mtd}]\label{prop:RiCor}
  Let $k$ be a commutative ring, and let $\sD_k$ denote the bicategory of
  $k$-algebras and derived categories of bimodules.  Suppose $B$ is a
  differential graded $k$-algebra, and let $T$ be a differential graded
  right $B$-module.  Let $A = \End_B(T)$ be the endomorphism algebra of
  $T.$ If $T$ has the following two properties, then $\sD_k(B)$ and
  $\sD_k(A)$ are equivalent as triangulated categories.
  \begin{enumerate}
  \item $T$ is a right-dualizable $B$-module.
  \item $T$ generates the triangulated category $\sD_k(B).$
  \end{enumerate}
\end{prop}
\begin{note}
  Since $T$ is right-dualizable, $T$ is a bounded complex of
  finitely-generated and projective $B$-modules.  Therefore the
  derived and underived endomorphism algebras are equal.
\end{note}
\begin{prop}[\cite{SS03Stable}]\label{prop:RiSpectra}
  Let $k$ be a commutative ring spectrum, and let $\sD_k$ denote the
  bicategory of $k$-algebras and homotopy categories of bimodules. Suppose
  $B$ is a $k$-algebra, and let $T$ be a right $B$-module.  Let
  $A = F_B(T,T)$ be the endomorphism ring spectrum of $T.$ If $T$ has the
  following two properties, then $\sD_k(B)$ and $\sD_k(A)$ are equivalent
  categories.
  \begin{enumerate}
  \item $T$ is right-dualizable as a $B$-module.
  \item $T$ generates the triangulated category $\sD_k(B).$
  \end{enumerate}
\end{prop}

\subsection{Acknowledgements} 
This work developed from a portion of the author's Ph.D.\ thesis.  The
author thanks his advisor, Peter May, for reading it many times in
draft form and offering a number of suggestions both mathematical
and stylistic.  The author also thanks Michael Ching for his
significant encouragement.  

The exposition has benefited from the helpful suggestions of anonymous
referees, and from the audiences who have listened to talks and read
drafts of this article.  The terminology ``Eilenberg-Watts
equivalence'' was suggested by Justin Noel during one of many helpful
conversations.

\section{Motivation}\label{sec:Mot-Exa}

To help motivate the generality of \fref{sec:context,sec:trbicat}
needed for \cref{thm:AzumayaEquiv,thm:Azumaya-char-triang}, we give
parallel introductions of the classification theorems for classical algebra and for
homotopical settings (\fref{sec:classical-brauer-gps,sec:Der-Cats},
respectively).  Both are special cases of the general theory.  In
\fref{eg:symm-mon-cats} we sketch a number of additional examples to
which the theory applies.

\subsection{Brauer groups of discrete rings}\label{sec:classical-brauer-gps}

In this section we review the classical theory of Brauer groups for
commutative (discrete) rings.  Let $k$ be a commutative ring and let
$\sB = \sM_k$ be the bicategory of $k$-algebras and their bimodules.
A 0-cell $A$ of this bicategory is a $k$-algebra, and $A^e$ is the
enveloping $k$-algebra $A \otimes_k A^{\op}$.  Note that throughout this
section we regard $A$ as a right module over $A^e$ and a left module
over $k$.
%
%
%
We recall the classical characterization/definition of Azumaya
algebras after recalling basic terminology and a well-known lemma.
These are special cases of the more general terms introduced in
\fref{defn:general-central-endomorphic,defn:general-fproj-fsep}.

\begin{samepage}
\begin{defns} \ 
\begin{itemize}
\item $A$ is called \emph{separable} over $k$ if $A$ is projective as a module
  over $A^e$.  Since $A$ is always finitely generated over $A^e$, this is
  equivalent to the condition that $A$ be a dualizable module over $A^e$.
  By the dual basis lemma, this is equivalent to the condition that the
  coevaluation
  \[
  A \otimes_{A^e} \Hom_{A^e}(A,A^e) \to \Hom_{A^e}(A,A)
  \]
  be an isomorphism. 
  To motivate this term, it should be noted that when $k$ is a field, $A$
  being separable over $k$ implies that $A$ is semi-simple over $k$ and
  remains semi-simple upon extension of scalars over any field extension of
  $k$.  When $A$ is also a field, this implies that $A$ is a separable
  extension of $k$ in the usual sense for fields \cite{Coh03Further}.
\item $A$ is called \emph{central} over $k$ if the center of $A$ is
  precisely $k$; this occurs if and only if the unit
  \[
  k \to \Hom_{A^e}(A,A)
  \]
  is an isomorphism.
\item $A$ is called \emph{faithfully projective} over $k$ if both the
  coevaluation
  \[
  \Hom_k(A,k) \otimes A \to \Hom_k(A,A)
  \]
  and the evaluation
  \[
  A \otimes_{A^e} \Hom_k(A,k) \to k
  \]
  are isomorphisms.  Note, again by the dual basis lemma, that the
  coevaluation being an isomorphism is equivalent to $A$ being
  finitely-generated and projective as a $k$-module.  The evaluation
  map being an isomorphism implies that $- \otimes_{k} A$ is
  \emph{object faithful}, meaning that $M \otimes_{k} A = 0$
  implies $M = 0$.
\end{itemize}
\end{defns}
\end{samepage}

\begin{lem}[See e.g.~\cite{KO74Theorie} or \cite{DI71Separable}]
  If $A$ is central and separable over $k$, then the evaluation
  \[
  \Hom_{A^e} (A, A^e) \otimes_k A \to A^e
  \]
  is an isomorphism.
\end{lem}

Here is the classical theorem, which follows directly from
\cref{thm:AzumayaEquiv}.  It provides a link between Eilenberg-Watts
equivalences and the more calculable conditions of central, separable and
faithfully projective.

\begin{thm*}[See \ref{thm:AzumayaEquiv}]\label{thm:characterization-thm-classical}
  Let $k$ be a commutative ring, and $A$ a $k$-algebra.  The following
  are equivalent:
  \begin{enumerate}
  \item $A$ is invertible as a $(k, A^e)$-bimodule, thus providing an
    Eilenberg-Watts equivalence between $A^e$ and $k$.
  \item $A$ is central and separable over $k$.
  \item \label{AzumayaEquivRtEC} $A$ is faithfully projective over $k$
    and the unit map
    \[
    A^e \to \Hom_{k} (A, A)
    \]
    is an isomorphism.
  \item \label{AzumayaEquivLtEC} There is a $k$-algebra $B$ such that $A \otimes_k B$ is Eilenberg-Watts
    equivalent to $k$.
  \end{enumerate}
\end{thm*}

\begin{rmk}
  Note that the last condition of \cref{AzumayaEquivRtEC} means that
  every $k$-linear endomorphism of $A$ is given by left and right
  multiplication in $A$.  We call such an algebra ``endomorphic'' in
  \cref{defn:central-endomorphic}.
\end{rmk}

\subsection{Homotopical Brauer groups}
\label{sec:Der-Cats}

In \fref{sec:hty-bicat} we describe a bicategorical structure on the
collection of categories $\sD_k(A,B)$, where $k$ is a commutative
($E_\infty$) d.g.~algebra or ring spectrum, cofibrant as a module over
itself, and $A$ and $B$ are cofibrant as $k$-modules (see
\fref{prop:homotopy-bicat-applications}).  This consists mainly of
assembling the necessary results from the literature on monoidal model
categories, but we have presented a careful treatment of the
Eilenberg-Watts structure.

In this section we continue to motivate the general theory by
describing the definitions and Azumaya characterization theorem in
this special case.  We use $\sma$ to denote the derived tensor
product, and $F$ to denote the derived hom.  The isomorphisms in the
homotopy category are the weak equivalences, and so we switch to this
terminology.  Let $k$ be a commutative ring, and let $A$ be a
differential-graded algebra over $k$.  For the topologically inclined
reader, $k$ and $A$ can be taken to be structured ring spectra---the
modern foundations of spectra such as \cite{EKMM97,MMSS01Model} allow
us to treat the algebraic and topological cases with the same
categorical arguments.  For either the algebraic or topological cases,
we assume throughout that $k$ and $A$ are cofibrant with respect to a
given model structure on $\sC_k$.

To begin, we make the following definitions inspired by the classical
case.  They are again special cases of the general bicategorical
definitions in \ref{defn:general-central-endomorphic} and
\ref{defn:general-fproj-fsep}.

\begin{defn}[Central/Endomorphic]
  \label{defn:central-endomorphic}
  We say that $A$ is \emph{central} if the natural unit map
  \[
  k \to F_{A^e}(A,A)
  \]
  is a weak equivalence.  Note that $F_{A^e}(A,A)$ is the (topological)
  Hochschild cohomology of $A$ over $k$.

  We say that $A$ is \emph{endomorphic} if the natural unit map from
  $A^e$ to the endomorphism algebra
  \[
  A^e \to F_k(A, A)
  \]
  is a weak equivalence.
\end{defn}

\begin{defn}[Faithfully Projective/Faithfully Separable]
  \label{defn:ffproj-ffsep}
  We say that $A$ is \emph{faithfully projective}
  over $k$ if both the
  coevaluation
  \[
  F_k(A, k) \sma_k A \to F_k(A, A)
  \]
  and the evaluation
  \[
  A \sma_{A^e} F_k(A, k) \to k
  \]
  are weak equivalences.  Note that the first of these two conditions is
  equivalent to requiring that $A$ be dualizable (in the derived category)
  over $k$.  (Duality is discussed in \cref{sec:dltybicat}.)  The second of
  these two conditions implies a kind of faithfulness over $k$, as in the
  classical case.

  We say that $A$ is \emph{faithfully separable} over $k$ if both the
  coevaluation
  \[
  A \sma_{A^e} F_{A^e}(A, A^e) \to F_{A^e}(A, A)
  \]
  and the evaluation
  \[
  F_{A^e}(A, A^e) \sma_k A \to A^e
  \]
  are weak equivalences.  The first of these two conditions is equivalent
  to requiring that $A$ be dualizable over $A^e$, and the second is a
  kind of faithfulness over $A^e$.
\end{defn}

\begin{defn}
  We call $A$ an \emph{Azumaya object} if any of the equivalent conditions
  below hold.
\end{defn}

\begin{samepage}
\begin{thm*}[See \ref{thm:AzumayaEquiv}]
 The following statements are equivalent:
  \begin{enumerate}
  \item $A$ is an invertible bimodule.
  \item $A$ is central and faithfully separable over $k$.
  \item \label{Azumaya-char-ffproj-sand}
    $A$ is endomorphic and faithfully projective over $k$.
  \item $A^e$ is Eilenberg-Watts equivalent to $k$.
  \item There exists a $k$-algebra $B$ such that
    $A \sma_k B$ is Eilenberg-Watts equivalent to $k$.
\end{enumerate}
\end{thm*}
\end{samepage}

\subsection{Additional Examples}

For concreteness, we focus our main applications on the
Eilenberg-Watts bicategories coming from differential graded algebras
or from ring spectra.  But the bicategorical setting is quite general
and here we sketch several additional examples to which the theory
applies.

\begin{example}[Symmetric monoidal categories]
\label{eg:symm-mon-cats}

The example of rings and modules motivates a standard construction for
a bicategory $\sM_{\sC}$ from any complete and cocomplete closed
symmetric monoidal category $\sC$.  The 0-cells are taken to be the
monoids of $\sC$, and 1-cells are taken to be bimodules.  Bimodule
maps are defined just as in the case of bimodules over rings, and
these constitute the 2-cells of $\sM_{\sC}$.  Colimits from $\sC$ are
used to construct the horizontal composition in $\sM_{\sC}$, in the
same way that coequalizers give the tensor product of two modules over
a ring.  Limits give $\sM_{\sC}$ a closed structure in the same way
that equalizers give the homomorphisms of two modules over a ring.

Remembering the underlying symmetric monoidal product of $\sC$ gives
$\sM_{\sC}$ an autonomous symmetric monoidal structure, and
$\sM_{\sC}$ is therefore Eilenberg-Watts.  \cref{thm:AzumayaEquiv}
thus characterizes the Azumaya monoids of $\sC$, and defines a Brauer
group of $\sC$ which reduces to the Brauer group of $k$ when $\sC$ is
the symmetric monoidal category of modules over a commutative ring
$k$.

The Brauer group of a symmetric monoidal category has been defined directly
by Pareigis \cite{Par75Non-Additive} and, for a braided monoidal category,
by Van Oystaeyen--Zhang \cite{VOZ98Brauer}.  These definitions (in the
symmetric monoidal case) are equivalent to ours, but there the emphasis is
on other aspects of the classical theory (e.g.\ separability and Casimir
elements); their work does not give the full characterization of Azumaya
algebras presented in \cref{thm:AzumayaEquiv}.  Likewise, we do not
discuss the variants on Brauer theory described by Pareigis or Van
Oystaeyen--Zhang.  Vitale \cite{Vit96Brauer} gives a thorough
treatment of Brauer groups for symmetric monoidal categories, and a
pleasant summary of categorical approaches to this case.
\end{example}

\begin{example}[Ex-spaces or parametrized spectra]
  May--Sigurdsson introduce $\Ex$, the bicategory of parametrized
  spectra in \cite[17.1.3]{may2006pht}.  The 0-cells of $\Ex$ are
  topological spaces, and the category of 1-cells $\Ex (A,B) = Ho
  \sS_{A \times B}$ is the homotopy category of \emph{parametrized
    spectra} over $A \times B$.  They describe a natural horizontal
  composition
  \[
  \odot : \Ex(A,B) \times \Ex(B,C) \to \Ex(A,C)
  \]
  by pullback along the diagonal $A \times B \times C \to A \times B
  \times B \times C$ and pushforward along the projection $A \times B
  \times C \to A \times C$.  These have adjoints given by pushforward
  and pullback along the same maps, and this defines a closed
  structure adjoint to the horizontal composition.

  The external smash product
  \[
  \bar{\sma} : \text{Ho} \sS_A \times \text{Ho} \sS_{A'} \to \text{Ho} \sS_{A \times A'}
  \]
  gives $\Ex$ a symmetric monoidal structure.  The ``opposite'' of a
  space $A$ is again $A$, and the equivalence $\Ex(A,B) \hty
  \Ex(B^{\op},A^{\op}) = \Ex(B,A)$ is induced by the switch map
  \[
  t: A \times B \to B \times A.
  \]
  Thus $\Ex$ is an Eilenberg-Watts bicategory and the characterization
  of Azumaya objects in $\Ex$ follows as in the case of spectra.
\end{example}

\begin{example} [Rings with many objects]
  Let $\aA$ be a small category enriched in $\Ab$, the category of
  abelian groups.  A left $\aA$-module is a functor $\aA \to \Ab$, and
  a right module is a functor $\aA^{\op} \to \Ab$; morphisms of
  $\aA$-modules are enriched natural transformations.  Now let $\Roid$
  be the bicategory whose 0-cells are small $\Ab$-categories, 1-cells
  are bimodules, and 2-cells are natural transformations.  The tensor
  product and hom for bimodules induce horizontal composition with a
  closed structure, and the tensor product in $\Ab$ induces a
  symmetric monoidal product on $\Roid$ making it an Eilenberg-Watts
  bicategory.  The unit object is the one-object category whose
  endomorphism ring is the ring of integers.

  Clearly this construction may be generalized by replacing $\Ab$ with
  a bicomplete closed symmetric monoidal category $\V$.  We denote the
  bicategory so constructed by $\V\text{--}\Roid$.  This too is an
  Eilenberg-Watts bicategory, and thus \fref{thm:AzumayaEquiv} gives a
  characterization of Azumaya objects.

  In the case that $\V = Mod_k$ is the category of modules over a commutative
  ring $k$, there is a natural morphism of bicategories
  \[
  \iota: \sM_k \to Mod_k\text{--}\Roid
  \]
  given by sending a $k$-algebra $A$ to the corresponding 1-object
  $Mod_k$-category, sending an $(A,B)$-bimodule to the functor whose
  value on the single object of $\iota(A) \times \iota(B)^{\op}$ is
  the given bimodule, and sending a morphism of bimodules to the
  induced natural transformation.  The unit object is $\iota(k)$, the
  one-object category with endomorphism ring $k$.  The morphism
  $\iota$ preserves the Azumaya property, and induces a homomorphism
  of Brauer groups
  \[
  Br(\sM_k) \to Br(Mod_k\text{--}\Roid).
  \]
For each pair of $k$-algebras, $A$ and $B$, the functor $\iota:
\sM_k(A,B) \to Mod_k\text{--}\Roid (\iota A , \iota B)$ is an
equivalence of categories and hence the induced map on Brauer groups
is an injection.  This map is in fact an isomorphism, as shown by
\cite[5.7]{borceux2002ac}.
\end{example}

\begin{example}[Sheaves]
The Brauer group of a sheaf $\sO$ of commutative rings was introduced
by Auslander in \cite{Aus66Brauer} and by Grothendieck in
\cite{grothendieck1968gba}.  By \fref{thm:AzumayaEquiv} this is the
same as the Brauer group of the symmetric monoidal category of
$\sO$-modules: a map of sheaves is an isomorphism if and only if it is
so locally, and thus an Azumaya object in the bicategory of
$\sO$-algebras and their bimodules is precisely the same as a sheaf of
Azumaya algebras.  This perspective was initiated by Auslander's work,
and has been discussed by Van Oystaeyen--Zhang \cite[3.5(2);
3.7]{VOZ98Brauer}.
\end{example}

\section{Eilenberg-Watts bicategories}
\label{sec:context}

Recall that we use the term ``Eilenberg-Watts bicategory'' for a
closed autonomous symmetric monoidal bicategory.  We describe these
terms in \fref{sec:monbicat,sec:symm-str,sec:auto-str,sec:Closed-Str}.
\fref{sec:dltybicat} reviews duality and invertibility in
bicategories.  In \fref{sec:AzumayaBrauer} we prove the
characterization theorem for Azumaya objects in Eilenberg-Watts
bicategories.

\begin{notn}
  We use arrows such as $f: M \to M'$ to denote that $f$ is a $2$-cell
  with source $M$ and target $M'$, and slashed arrows such as
  $\cell{M}{A}{B}$ to denote that $M$ is a $1$-cell with source $A$
  and target $B$.  We use $\circ$ or juxtaposition to denote vertical
  composition of $2$-cells, and $\odot$ to denote horizontal
  composition of $1$-cells and of $2$-cells.  We write
  $\odot$-composition in ``diagrammatic order'', so that a composite
  of 1-cells
  \[
  A \sfto{M} B \sfto{N} C
  \]
  is denoted $\cell{(M \odot N)}{A}{C}.$
\end{notn}

\subsection{Monoidal bicategories}\label{sec:monbicat} 
A monoidal bicategory can be defined as a tricategory with one object,
in the sense of \cite{GPS95Coherence}.  In practical terms, this means
that the bicategory is equipped with an additional monoidal product on
0-, 1-, and 2-cells forming a pseudofunctor
\[
  \otimes \cn \sB \times \sB \to \sB.
\]
For categories of 1-cells $A \sto B$ and $A' \sto B'$, this
gives functors
\[
  \otimes \cn \sB(A,B) \times \sB(A', B') \to \sB(A \otimes A', B
  \otimes B').
\]
For associativity and unit constraints we have invertible
1-cells
\[
\cell{a}{(A_1 \otimes A_2) \otimes A_3}{A_1 \otimes (A_2 \otimes A_3)}
\]
and
\begin{align*}
  \cell{l &}{A \otimes k}{A}\\
  \cell{r &}{k \otimes A}{A}.
\end{align*}
Invertible 2-cells fill the associativity pentagon and unit triangle,
and there are two additional invertible 2-cells relating the left and
right units with the associativity constraint.  This data satisfies a
number of reasonable axioms which we will not need to reference
explicitly.  

In lieu of further details, which can be found in
\cite{GPS95Coherence,Str03Functorial,Sch09Classification}, we
encourage the reader unfamiliar with monoidal bicategories to
concentrate on the representative example of algebras and bimodules.
Here the 0-cells are algebras over a commutative ring $k$, 1-cells
$A \sto B$ are bimodules, and 2-cells are bimodule homomorphisms.
The horizontal composition
\[
A \sto B \sto C
\]
is given by the tensor product of bimodules $\otimes_B$.  The tensor
product of algebras $\otimes_k$ gives a monoidal structure on this
bicategory.  This example contains all of the essential structure
present in the more sophisticated algebraic and topological examples,
and also serves to illustrate the additional bicategorical structure
described below.

A trivial but slightly technical lemma about monoidal
structure will be useful later:
\begin{lem}\label{lem:bicat-unit-comp}
  Let $k$ be the unit of a monoidal bicategory, and suppose that
  $\cell{M}{C}{k}$ and $\cell{N}{k}{D}$ are 1-cells.  Then there are
  natural isomorphisms
  \begin{align*}
    N \otimes M \iso & \  M \odot N\\
    N \otimes M \iso & \ (N \otimes C) \odot (D \otimes M)
  \end{align*}
  and hence there is a 2-cell isomorphism filling the square below:
  \[
  \cmxymat{
    C \ar[r]^-{N \otimes C}|-@{|} \ar[d]_-{M}|-@{|} & D \otimes C \ar[d]^-{D \otimes M}|-@{|}
    \\
    k \ar[r]^-{N}|-@{|} & D
  }
  \]
\end{lem}
\begin{proof}
  Pseudofunctorality of the monoidal product means we have isomorphisms
  \begin{align*}   
    (k \odot N) \otimes (M \odot k)
    \iso & \ 
    (k \otimes M) \odot (N \otimes k) 
    \\
    \intertext{and}
    (N \odot D) \otimes (C \odot M)
    \iso & \ 
    (N \otimes C) \odot (D \otimes M).
  \end{align*}
  Composing with unit isomorphisms then gives the result.
\end{proof}


\subsection{Symmetric monoidal bicategories}
\label{sec:symm-str}

Symmetry for the monoidal product in a bicategory entails invertible
1-cells
\[
\cell{R_{AB}}{A \otimes B}{B \otimes A},
\]
and again the example of algebras and bimodules is illuminating.  The
structure 1-cells form a braiding and for each pair $(A,B)$ there is
an invertible 2-cell (syllepsis) making the composite 1-cell
\[
R_{AB} \odot R_{BA} \cn A \otimes B \sto B \otimes A \sto A \otimes B
\]
isomorphic to the identity 1-cell.  The symmetry for this syllepsis is
an axiom equating the two isomorphisms $R_{AB} \odot
R_{BA} \odot R_{AB} \iso R_{AB}$, one coming from the syllepsis for
$(A,B)$ and the other coming from that of $(B,A)$.  Descriptions of
braided, sylleptic, and symmetric monoidal bicategories, together with
coherence theorems in the braided and symmetric cases, are found in
\cite{Gur11Loop,GO12Infinite}.

\subsection{Autonomous structure}\label{sec:auto-str}
A symmetric monoidal bicategory $\sB$ is \emph{autonomous} \cite{DaS1997}
if each 0-cell $A$ has an opposite, $A^{\op}$, providing an involution on
$\sB$ and such that there are (suitably natural) equivalences of categories
\[
\sB(A \otimes B,C) \hty \sB(B,A^{\op} \otimes C) \text{\quad  and \quad}
\sB(B \otimes A^{\op}, C) \hty \sB(B, C \otimes A)
\] 
for all $B$ and $C$.  As in classical ring theory, where $A^{\op}$ is the
opposite $k$-algebra, we write $A^e$ for the enveloping object
$A \otimes A^{\op}$, and given a 1-cell $\cell{M}{A}{B}$ we write
$\cell{M^{\op}}{B^{\op}}{A^{\op}}$ for its image under the equivalence
\[
\sB(A,B) \hty \sB(B^{\op},A^{\op}).
\]

\subsection{Closed structure}\label{sec:Closed-Str}
A \emph{closed structure} for a bicategory, $\sB$, defines right adjoints
for $\odot$.  For a $1$-cell $M$, the right adjoint to $M \odot -$ is
called ``right-hom'', or ``source-hom'', and denoted $\shom{M}{-}$.  The
adjoint to $- \odot M$ is called ``left-hom'', or ``target-hom'', and
denoted $\thom{M}{-}$.  For $M \in \sB(A,B)$ and $W \in \sB(A,C),$ we
have the 1-cell $\shom{M}{W} \in \sB(B,C).$ For $U \in \sB(D,B),$ we have
the 1-cell $\thom{M}{U} \in \sB(D,A).$ The adjunctions are written as
\[\begin{array}{ccc}
  \sB(M \odot V,W) & \iso & \sB(V,\shom{M}{W})\\
  \sB(T \odot M,U) & \iso & \sB(T,\thom{M}{U})
\end{array}\]

The existence of left and right hom functors defines a \emph{closed
  bicategory}.  A thorough description of closed structures can be
found in \cite{may2006pht}.


\subsection{Duality and invertibility in bicategories}\label{sec:dltybicat} 
For general discussion about duality, we consider fixed 1-cells
$\cell{X}{A}{B}$ and $\cell{Y}{B}{A}$ in a closed bicategory $\sB$.

\begin{defn}[Dual pair]
  We say $(X,Y)$ is a dual pair, or ``$X$ is left-dual to $Y$'' (``$Y$ is
  right-dual to $X$''), or ``$X$ is right-dualizable'' (``$Y$ is
  left-dualizable'') to mean that we have 2-cells
  \[
  \eta :A \to X \odot Y \text{\quad and \quad} \epz: Y \odot X \to B
  \]
  such that the following composites are the respective identity
  2-cells.  These are known as the triangle identities:
  \[
  X \iso A \odot X \fto{\eta \odot \id} X \odot Y \odot X \fto{\id
    \odot \epz} X \odot B \iso X
  \]
  \[
  Y \iso Y \odot A \fto{\id \odot \eta} Y \odot X \odot Y \fto{\epz
    \odot \id} B \odot Y \iso Y.
  \]
\end{defn}

\begin{lem}\label{lem:duality-composition}
  Let $k$ be the unit of a symmetric monoidal category, and $A$ a
  $k$-algebra. Let ${_k}A$, $A_k$ denote the base-change 1-cells
  induced by the unit map $k \to A$.  Then $({_k}A, A_k)$ is a dual
  pair with structure maps
  \[
  k \to {_k}A \otimes_A A_k
 \]
 and 
 \[
 A_k \otimes_k {_k}A \to A.
 \]
\end{lem}

\begin{note}
  When $(X,Y)$ is a dual pair, we will occasionally refer to this by saying
  that $X$ is dualizable \emph{over $B$}, since the unit condition often
  amounts to a finiteness of $X$ over $B$.  When $(Y,X)$ is a dual pair, we
  will say that $X$ is dualizable \emph{over $A$}.  With this convention,
  the phrase ``dualizable over'' always references the target of the
  evaluation map, and uniquely determines whether we mean left- or
  right-dualizable.
\end{note}

\begin{defn}[Invertible pair]\label{defn:invDefn}
  A dual pair $(X,Y)$ is called invertible if the maps $\eta$ and
  $\epz$ are isomorphisms.  Equivalently, the adjoint pairs
  described in \cref{prop:dltyAdj} are adjoint equivalences.
\end{defn}

Duality for monoidal categories has been studied at length; one
reference for duality in a bicategorical context is \cite[\S
16.4]{may2006pht}, and there are surely others in the categorical
literature.  The definition of duality does not require $\sB$ to be
closed, but we will make use of the following basic facts about
duality, some of which do require a closed structure on $\sB$.  The
following two results can be found in \cite[\S 16.4]{may2006pht}.

\begin{prop}
  A 1-cell $X \in \sB(A,B)$ is right-dualizable if and only if the coevaluation
  \[
  \nu: X \odot \thom{X}{B} \to \thom{X}{X}
  \]
  is an isomorphism.  Moreover, this is the case if and only if the
  map
  \[
  \nu_Z: X \odot \thom{X}{Z} \to \thom{X}{X \odot Z}
  \] 
  is an isomorphism for all 1-cells $\cell{Z}{B}{B}.$
\end{prop}

\begin{samepage}
\begin{prop}\label{prop:dltyAdj} Let $(X,Y)$ be a dual pair in
  $\sB$, with $\cell{X}{A}{B}$ and $\cell{Y}{B}{A}.$
\begin{enumerate}
\item For any 0-cell $C$, we have two adjoint pairs of functors, with
  left adjoints written on top:
  \[\cmxymat{
    {\sB}(C,A) \ar[r]<.4ex>^-{- \odot X} & {\sB}(C,B) \ar[l]<.4ex>^-{- \odot Y}
  }\]
  \[\cmxymat{
    {\sB}(A,C) \ar[r]<.4ex>^-{Y \odot -} & {\sB}(B,C) \ar[l]<.4ex>^-{X \odot -}
  }\]
  The structure maps for the dual pair give the triangle identities
  necessary to show that the displayed functors are adjoint pairs.

\item The right dual, $Y,$ is canonically isomorphic to $\thom{X}{B}$, and
  for any 1-cell $\cell{W}{D}{B}$, the coevaluation map
  $W \odot \thom{X}{B} \to \thom{X}{W}$ is an isomorphism.

\item The left dual, $X,$ is canonically isomorphic to $\shom{Y}{B}$, and
  for any 1-cell $\cell{U}{A}{D'}$, the coevaluation map
  $\shom{Y}{A} \odot U \to \shom{Y}{U}$ is an isomorphism.
\end{enumerate}
\end{prop}
\end{samepage}

\begin{lem}\label{lem:dlty-eval}
  Let $\cell{X}{A}{B}$ be a 1-cell in $\sB(A,B)$.  If $X$ is
  right-dualizable and the unit $A \to \thom{X}{X}$ is an isomorphism,
  then the evaluation $X \odot \shom{X}{A} \to A$ is an isomorphism.
  Likewise, if $X$ is left-dualizable and the unit $B \to \shom{X}{X}$
  is an isomorphism, then the evaluation $\thom{X}{B} \odot X \to B$
  is an isomorphism.
\end{lem}
\begin{proof}
  We prove the first statement, leaving the second as an exercise
  in opposites.  Let $Y$ denote the right dual of $X$.  Since $X$
  is right-dualizable, $Y$ is left-dualizable and $X$ is
  isomorphic to the canonical left dual of $Y$:
  $X \iso \shom{Y}{B}$.  The isomorphism
  $A \fto{\iso} \thom{X}{X}$ implies that the unit for the duality
  is an isomorphism: $A \fto{\iso} X \odot Y$.  Now we have the
  following commutative square:
  \[
  \cmxymat{
    X \odot_B \shom{X}{A} \ar[rr]^-{\text{evaluation}} \ar[d]_-{\iso} 
    & 
    & A \ar[d]^-{\iso} \\
    \shom{Y}{B} \odot_B \shom{X}{A} \ar[r]^-{\iso} 
    & \shom{Y}{\shom{X}{A}} \ar[r]^-{\iso} 
    & \shom{(X \odot Y)}{A}
  }
  \]
  where the two vertical isomorphisms are described above, the
  left-hand isomorphism is a consequence of dualizability for $Y$, and
  the right-hand isomorphism is an exercise in adjunction.
\end{proof}

\subsection{Characterization of Azumaya objects} 
\label{sec:AzumayaBrauer}
In this section we give a proof of \cref{thm:AzumayaEquiv} in two
stages.  Let $A$ be a fixed 0-cell of an Eilenberg-Watts bicategory
$\sB$.

\begin{notn}
  Throughout, it will be useful to distinguish between the 0-cell $A$,
  the unit 1-cell $\cell{A}{A}{A}$ and the 1-cell $\cell{A}{k}{A^e}$
  corresponding to $\cell{A}{A}{A}$ under the autonomy $\sB(A,A) \hty
  \sB(k,A^e)$.  We will refer to the latter as the \emph{floor} of
  $A$, and denote it as $\fl A$, with the motivation that $\fl A$ is
  naturally induced by $A$, but with its source being the ground ring
  and its target being the enveloping object, $A^e$.  Thus $\fl A$ is a 1-cell
  with source $k$ and target $A^e$:
  \[
  \cell{\fl A}{k}{A^e}.
  \]
  We will not make explicit reference to the related 1-cell $A^e \sto
  k$, but a compatible notation would be the \emph{ceiling}, $\lceil A
  \rceil$.
  As for the 0-cell $A$ and corresponding unit 1-cell, these sit
  at different levels of the ambient bicategory, and this will be
  sufficient to distinguish them in context.
\end{notn}

We now give the general characterization of Azumaya objects;
\cref{thm:AzumayaEquiv} follows directly from this.  First, we give
natural generalizations of the classical terminology:

\begin{defn}[Central/Endomorphic]\label{defn:general-central-endomorphic}
  We say that $A$ is \emph{central} over $k$ if the natural unit map
  \[
  k \to \thom{\fl A}{\fl A}
  \]
  is an isomorphism.

  We say that $A$ is \emph{endomorphic} over $k$ if the natural unit map
  \[
  A^e \to \shom{\fl A}{\fl A}
  \]
  is an isomorphism.
\end{defn}

\begin{defn}[Faithfully Projective/Faithfully Separable]
  \label{defn:general-fproj-fsep}
  We say that $A$ is \emph{faithfully projective}
  over $k$ if both the
  coevaluation
  \[
  \shom{\fl A}{k} \sma_k \fl A \to \shom{\fl A}{\fl A}
  \]
  and the evaluation
  \[
  \fl A \sma_{A^e} \shom{\fl A}{k} \to k
  \]
  are isomorphisms.

  We say that $A$ is \emph{faithfully separable} over $k$ if both the
  coevaluation
  \[
  \fl A \sma_{A^e} \thom{\fl A}{A^e} \to \thom{\fl A}{\fl A}
  \]
  and the evaluation
  \[
  \thom{\fl A}{A^e} \sma_k \fl A \to A^e
  \]
  are isomorphisms.
\end{defn}

We are now ready to prove the characterization theorem for Azumaya
objects in Eilenberg-Watts bicategories.

\begin{proof}[Proof of \cref{thm:AzumayaEquiv}]
  In general, a pair of 1-cells $(X,Y)$ is invertible if and only if $X$ is
  right-dualizable, $Y$ is isomorphic to the canonical right dual of $X$,
  and both the unit and counit of the duality are isomorphisms.  This gives
  the equivalence of \ref{AzumayaEquivInv} and
  \ref{AzumayaEquivRtEC}. Likewise, $(Y,X)$ is invertible if and only if
  $X$ is left-dualizable, $Y$ is isomorphic to the canonical left dual of
  $X$, and both the unit and counit of the duality are isomorphisms.  This
  gives the equivalence of \ref{AzumayaEquivInv} and
  \ref{AzumayaEquivLtEC}.  Since $(X,Y)$ is an invertible pair if and only
  if $(Y,X)$ is such, the first three conditions are seen to be equivalent.

  Connecting the last two conditions to the previous ones is
  essentially a bicategorical diagram chase: of course if $\fl A$ is
  invertible, then the implication of the last two conditions is
  clear.  For the reverse, suppose that $\cell{P}{k}{A \otimes B}$ is
  an invertible 1-cell between $k$ and $A \otimes B$ and let $P^*$ be
  its inverse.  Let $P^e = P \otimes (P^*)^{\op}$ and $P^{*e} = P^*
  \otimes P^{\op}$.  Then $(P^e,P^{*e})$ is an invertible pair of
  1-cells between the 0-cells $k^e$ and $(A \otimes B)^e \iso A^e
  \otimes B^e$ (using the symmetry).  Moreover, since $(P,P^*)$ is an
  invertible pair we have $P^e = P \otimes (P^*)^{\op} \iso \fl{A \otimes B}
  \iso \fl A \otimes \fl B$ (again using symmetry).  The diagram below
  describes our situation:
  \[
  \cmxymat@R=3pc@C=3pc{ \; \phantom{k}k \;
    \ar@/^.4pc/[rr]<.4ex>^-{P^e}|-{\object@{|}}
    \ar@/^2.5pc/[rr]^-{\fl A \otimes \fl B}|-{\object@{|}}
    \ar@/_/[dr]_-{\fl A}|-{\object@{|}} & &
    (A \otimes B)^e \ar@/^.4pc/[ll]<.4ex>^-{P^{*e} = Q}|-{\object@{|}}\\
    & A^e \ar@/_/[ur]_-{A^e \otimes \fl B}|-{\object@{|}} & }\]

  Let $Q = P^{*e}$ and let $T = (A^e \otimes \fl B) \odot Q$.  We have
  shown that $(\fl A \otimes \fl B , Q)$ is an invertible pair of
  1-cells, and we will now show that $(\fl A, T)$ is an invertible
  pair between $A^e$ and $k$.  A straightforward use of associativity
  shows that $\fl A \odot T \iso k$.  We now use
  \fref{lem:bicat-unit-comp} for the reverse composite.  This shows
  that $\fl A$ is invertible, with the triangle identities coming from
  those of $(P^e, P^{*e})$:
\begin{align*}
    T \odot \fl A =(A^e \otimes \fl B) \odot Q \odot \fl A 
    & \iso
    (A^e \otimes \fl B) \odot (\fl A \otimes A^e \otimes B^e) \odot (A^e \otimes Q)\\
    & \iso 
    (A^e \otimes \fl B) \odot (A^e \otimes \fl A \otimes B^e) \odot (A^e \otimes Q)\\
    & \iso 
    (A^e \otimes \fl A \otimes \fl B) \odot (A^e \otimes Q)\\
    & \iso (A^e \otimes k) \iso A^e.
  \end{align*}
\end{proof}

\section{Triangulated bicategories}\label{sec:trbicat}

We recall first the definitions of localizing subcategory and
generator for a triangulated category, and then give a definition
(\ref{defn:trbicatDef}) of triangulated bicategory suitable for our
purposes.  In particular, under this definition $\sD_k$ is a
triangulated bicategory when $k$ is a commutative d.g.~ring or ring
spectrum.

\begin{defn}[Localizing subcategory]
  If $\sT$ is a triangulated category with infinite coproducts, a
  \emph{localizing} subcategory, $\sS$, is a full triangulated
  subcategory of $\sT$ which is closed under coproducts from $\sT$.
\end{defn}

\begin{rmk}
  This is equivalent to the definition for arbitrary triangulated
  categories of \cite{hovey1999mc}, (which requires that a localizing
  subcategory be thick) because a triangulated subcategory
  automatically satisfies the 2-out-of-3 property and because in any
  triangulated category with countable coproducts, idempotents have
  splittings.  See \cite[1.5.2, 1.6.8, and 3.2.7]{neeman2001tc} for details.
\end{rmk}

\begin{defn}[Triangulated generator]\label{defn:trGen}
  A set, $\sP$, of objects in $\sT$ (triangulated category with
  infinite coproducts, as above) is a set of \emph{triangulated
    generators} (or simply \emph{generators}) if the only localizing
  subcategory containing $\sP$ is $\sT$ itself.
\end{defn}

\begin{defn}[Triangulated bicategory {\cite[\S 16.7]{may2006pht}}]  
  \label{defn:trbicatDef}
  A closed bicategory $\sB$ will be called a \emph{triangulated bicategory}
  if for each pair of 0-cells, $A$ and $B$, $\sB(A,B)$ is a triangulated
  category with infinite coproducts, and furthermore the local
  triangulations on $\sB$ are compatible as described in the following
  axioms.
\begin{itemize}
\item[(TC0)] The local suspension functors 
  \[
  \Sigma_{A,B} : \sB(A,B) \fto{\hty} \sB(A,B)
  \]
  assemble as the components of a strong transformation
  \[
  \Sigma : \sB(-,-) \to \sB(-,-).
  \]
  This implies, in particular, that for a 1-cell $\cell{X}{A}{B}$ there are
  natural structure isomorphisms
  \[\cmxymat{
    \Sigma A \odot X \ar[r]^-{\al_X}_-{\iso} 
    & \Sigma (A \odot X) \iso 
    \Sigma X \iso 
    \Sigma(X \odot B) 
    & X \odot \Sigma B. \ar[l]_-{\be_X}^-{\iso} 
  }\]
\item[(TC1)] For a 0-cell $A$, the following composite of structure
  isomorphisms interchanging suspension coordinates is multiplication by
  -1:
 \[
 \Sigma^2 A \fto{\Sigma \be_A^{-1}} \Sigma(A \odot \Sigma A) 
 \fto{\al_{\Sigma A}^{-1}}
 \Sigma A \odot \Sigma A 
 \fto{\be_{\Sigma A}} 
 \Sigma(\Sigma A \odot A)
 \fto{\Sigma \al_A}
 \Sigma^2 A.
  \]

\item[(TC2)] Each of $- \odot -$, $\shom{-}{-}$, and $\thom{-}{-}$ is exact
  in both variables.
\end{itemize}
\end{defn}

The axioms above generalize those of a triangulated monoidal category
although, because we do not need it, we do not include a version of
the braid axiom (TC3).

\begin{rmk}
  \label{rmk:k-alg-trbicat}
  For a commutative d.g.~algebra or ring spectrum $k$, each local
  category $\sD_k(A,B)$ is triangulated, with the suspension functor
  given by horizontal composition $\Si A \odot -$.  The monoidal
  category $\sD_k(k,k)$ is well-known to satisfy the compatibility
  conditions, and the same reasoning shows that $\sD_k$ satisfies
  these conditions in general.
\end{rmk}

If $\sB$ is a triangulated bicategory and $P$, $Q$ are 1-cells in
$\sB(A,B)$, we emphasize that $\sB$ is triangulated by writing the
abelian group of 2-cells $P \to Q$ as $\sB[P,Q]$ and by writing the
graded abelian group obtained by taking shifts of $Q$ as $\sB[P,Q]_*$.
To emphasize the source and target of $P$ and $Q$, we may also write
$\sB(A,B)[P,Q]_*$.

\begin{defn}[$\odot$-faithful 1-cells]\label{defn:odD}
  In any locally additive bicategory, $\sB$, a 1-cell $\cell{W}{A}{B}$
  is called \emph{source-faithful} if triviality for any 1-cell
  $\cell{Z}{C}{A}$ is detected by triviality of the composite $Z \odot
  W$.  That is, $\cell{Z}{C}{A}$ is zero if and only if $Z \odot W =
  0$.  A collection of 1-cells, $\sE$, in $\sB(A,B)$ is called
  \emph{jointly source-faithful} if the objects have this property
  jointly; that is, $Z = 0$ if and only if $Z \odot W = 0$ for all $W
  \in \sE$.  The term \emph{target-faithful} is defined similarly,
  considering $W \odot -$ instead of $- \odot W$.
\end{defn}

\begin{rmk}
  If $\sC$ is an additive monoidal category with monoidal product $\odot$,
  the unit object is both source- and target-faithful.  In an arbitrary
  locally additive bicategory $\sB$, if $A \ne B$ then $\sB(A,B)$ may not
  have a single object with this property.  In relevant examples, however,
  the collection of all 1-cells, ob$\sB(A,B)$, does have this property
  jointly.
\end{rmk}

\begin{lem}\label{lem:odDlem}
  Let $\sB$ be a triangulated bicategory, and let $\cell{P}{A}{B}$ be a
  generator for $\sB(A,B)$.  If the collection of all 1-cells,
  $\sB(A,B)$, is jointly source-faithful (resp.\ target-faithful), then $P$ is
  source-faithful (resp.\ target-faithful).
\end{lem}
\begin{proof}
  Consider the source-faithful case; the target-faithful case is similar.
  Given any 1-cell 
  \[
  \cell{Z}{C}{A}
  \]
  with $Z \odot P = 0$, let $\sS$ be the full subcategory of 1-cells,
  $\cell{W}{A}{B}$ for which $Z \odot W = 0$.  This is a localizing
  subcategory of $\sB(A,B)$, and by assumption $P \in \sS$, so
  $\sS = \sB(A,B)$, and hence $Z = 0$.
\end{proof}

\begin{rmk}\label{rmk:odDrmk}
  Since the functors $P \odot -$ are exact, the property of $P \odot
  -$ detecting trivial objects is equivalent to $P \odot - $ detecting
  isomorphisms (meaning that a 2-cell $f$ is an isomorphism if and
  only if $P \odot f$ is so).
\end{rmk}

\subsection{Homotopy Bicategories}\label{sec:hty-bicat}
In this section, let $\sC$ be a complete and cocomplete closed
symmetric monoidal model category.  Let $\sM$ be the closed bicategory
formed by $\sC$-monoids and their bimodules as in
\cref{eg:symm-mon-cats}.  We say that $\sM$ has a local model
structure if each category of 1- and 2-cells $\sM(A,B)$ is a model
category.  

We now describe further conditions for the model structure on $\sC$
which ensure that $\sM$ has a local model structure and that the
collection of homotopy categories assembles to an Eilenberg-Watts
bicategory.  The first of these, pushout products, implies that the
horizontal composition of 1-cells descends to a homotopy bicategory.
The second, unit replacement, implies that horizontal composition of
1-cells on the level of homotopy is unital.

\begin{defn}[Pushout products]
  Suppose $\sM$ has a local model structure.  We say $\sM$ has
  pushout products if, for any 0-cells $A$, $B$, and $C$ and for any
  cofibrations $f \cn U \to V$ in $\sM(A,C)$ and $g \cn W \to X$ in
  $\sM(C,B)$, the pushout product $f \pop g$ below is a cofibration in
  $\sM(A,B)$ which is acyclic if either $f$ or $g$ is acyclic.

  \[
  \cmxymat@C=.3cm@R=.3cm{
    U \sma_C W \ar[rr] \ar[dd] & \hspace{1ex} & V \sma_C W \ar[dd] &\\
    & \vspace{2ex}  & \vspace{1ex} & \\
    U \sma_C X \ar[rr] & \hspace{1ex} &  \bP \ar@{..>}[dr]^-{f \pop g} & \\
    & & & V \sma_C X
  }
  \]
\end{defn}
\begin{defn}[Unit replacement condition]
  Let $A$ be a 0-cell of $\sM$, and let $QA$ be a cofibrant
  replacement for the unit 1-cell $A$ in $\sM(A,A)$.  We say that the
  unit replacement condition holds for $A$ if, for any cofibrant $X
  \in \sM(A,C)$, the induced map
  \[
  QA \odot X \to A \odot X \iso X
  \]
  is a weak equivalence in $\sM(A,C)$.  Note that this condition is
  independent of the choice of cofibrant replacement $QA$, and is
  automatically satisfied if $A$ is cofibrant in $\sB(A,A)$.
\end{defn}

Next we recall one result which implies that $\sM$ has a local model
structure, and then give another showing that this local model
structure indeed descends to form a homotopy bicategory.

\begin{prop}[{\cite[4.1]{SS00Algebras}}]\label{prop:ss-model-str}
  Let $\sC = (\sC, \sma, k)$ be a cofibrantly generated closed
  monoidal model category which satisfies the monoid axiom and in
  which each object is small with respect to the whole category.  Let
  $\sM_\sC$ be the bicategory formed by $\sC$-monoids which are
  cofibrant in $\sC$ and their bimodules.  Then $\sM_\sC$ has a local
  model structure, and has pushout products.  Moreover, the category
  of $\sC$-monoids $\Mon(\sC)$ is a cofibrantly generated model
  category, and if $k$ is cofibrant then every cofibration in
  $\Mon(\sC)$ whose source is cofibrant is also a cofibration upon
  forgetting to $\sC$.
\end{prop}

This proposition applies, for example, to simplicial sets,
$\Ga$-spaces, symmetric spectra, simplicial abelian groups, chain
complexes, and $S$-modules.  Details for these and other examples are
given in \cite[\S5]{SS00Algebras}.

\begin{prop}[Homotopy bicategory]\label{prop:hty-bicat}
  Let $\sC$ be as above, and assume the following:
  \begin{itemize}
  \item The unit, $k$, is cofibrant in $\sC$.
  \item The unit replacement condition holds for each monoid in $\sC$
    which is cofibrant in $\sC$.
  \item For each monoid $A$ and cofibrant $A$-module $N$, $-
    \sma_A N$ takes weak equivalences of $A$-modules to weak
    equivalences in $\sC$.
  \end{itemize}
  Then the collection of homotopy categories
  $h\sM_\sC(A,B)$ forms a closed bicategory, which we denote $h\sM_\sC$.
\end{prop}
\begin{proof}
  We begin with an application of the theory of Quillen functors of
  two variables \cite[\S 4.2]{hovey1999mc}.  Since $\sM_\sC$ has pushout
  products by \fref{prop:ss-model-str}, each component of horizontal
  composition is a left Quillen bifunctor
  \[
  \odot \cn \sM(A,C) \times \sM(C,B) \to \sM(A,B)
  \]
  and thus induces a bifunctor on $h\sM$ which we also denote  with
  $\odot$.  The unit replacement condition ensures that $\odot$ is
  unital on $\sM_\sC$, and the third condition ensures that
  $h\sM_\sC(A,B) \hty h\sM_\sC(A',B')$ if $A \hty A'$ and $B \hty B'$
  are weak equivalences in $\Mon(\sC)$ \cite[Theorem 3.3]{SS00Algebras}.
  %
  %
\end{proof}

\begin{prop}\label{prop:model-bimod-bicat}
  Let $\sC$ be as in \fref{prop:hty-bicat}.  Then the derived monoidal
  product on $\sC$ descends to an autonomous symmetric monoidal
  structure on $h\sM$, making $h\sM$ Eilenberg-Watts.
\end{prop}
\begin{proof}
  If $A$ and $B$ are cofibrant in $\sC$, then $A \sma B$ is a
  model for the derived product, and a cofibrant $(A,B)$-bimodule is
  also cofibrant in $\sC$.  Moreover, the product in $\sC$ preserves
  bimodule structure and hence induces a bifunctor
  \[
  \sma: h\sM_\sC(A,B) \times h\sM_\sC(C,D) 
  \to h\sM_\sC(A \sma C, B \sma D)
  \]
  which agrees (for each $A, B, C, D \in \Mon(\sC)$) with the derived
  monoidal product in $\sC$.  This gives an associative product on the
  bicategory $h\sM_\sC$, which is unital because $k$ is cofibrant in
  $\sC$.  The autonomous structure on $h\sM_\sC$ is obtained by taking
  the opposite multiplication on monoids and modules of $\sC$, and the
  symmetry of the product in $\sC$ gives the symmetry in $h\sM_\sC$.
\end{proof}

We now give our main application, in which $k$ is taken to be a
commutative d.g.~algebra or commutative ring spectrum and $\sC =
\sC_k$ is taken to be the category of $k$-modules.  One may work with
any monoidal model category of spectra, but we must assume that $k$ is
cofibrant in $\sC_k$ and thus the result applies only to those
categories in which the unit is cofibrant (e.g., symmetric spectra or
orthogonal spectra).  Likewise, if $k$ is a d.g.~algebra, one must
choose a model structure on $\sC_k$ for which the unit is cofibrant
(e.g., the injective model structure).
\begin{prop}\label{prop:homotopy-bicat-applications}
  Let $k$ be a commutative d.g.~algebra or commutative ring spectrum,
  and let $\sC_k$ be the symmetric monoidal category of $k$-modules,
  with a model structure for which $k$ is cofibrant.  Then we have
  Eilenberg-Watts bicategories $\sM_k = \sM_{\sC_k}$ and $\sD_k =
  h\sM_{\sC_k}$.
\end{prop}
\begin{proof}
  The hypotheses of \fref{prop:ss-model-str,prop:hty-bicat} are
  satisfied by $\sC_k$, so \fref{prop:model-bimod-bicat} applies: most
  of the hypotheses are discussed in \cite{SS00Algebras} and so we
  comment only on the unit replacement condition.  Suppose $A \in
  Mon_c(\sC)$, and consider a specific cofibrant replacement for $A$
  as an $(A,A)$-bimodule: the two-sided bar construction $QA =
  B(A,A,A)$.  Then $QA \sma_A X \hty B(A,A,X)$, and the latter
  retracts to $A \sma_A X$ along the simplicial retraction induced by
  the unit of $A$.  Thus the induced map
  \[
  QA \sma_A X \to A \sma_A X \hty X
  \]
  is an equivalence in $\sD_k(A,A)$.
\end{proof}

\subsection{Invertibility in triangulated bicategories} 
\label{sec:invertibility-trbicat}
For this section, we let $\sD$ denote a triangulated Eilenberg-Watts
bicategory with unit 0-cell $k$.  Bousfield localization for
triangulated categories \cite[\S4.10]{Kra09Localization} generalizes
to triangulated bicategories in two ways, using 1-cell composition
over either the source or target of a given 1-cell.  Recall that we
write 1-cell composition in diagrammatic order:
\[
\left( C \sfto{M} A \right)
\odot
\left( A \sfto{T} B \right)
=
\left(
\cmxymat@R=3pc@C=3pc{
  C \ar[r]^-{M \odot T}|-{{\rule[0pt]{.2pt}{3pt}}}
  & B
}
\right)
\]

\begin{defn}
  Let $\cell{T}{A}{B}$ be a 1-cell in $\sD(A,B)$. 

  A 1-cell $\cell{M}{C}{A}$ is \emph{source-$T$-acyclic} if $M \odot T =
  0$.  A 1-cell $\cell{N}{C}{A}$ is \emph{source-$T$-local} if
  \[
  \sD(C,A)[M,N]_* = 0
  \]
  for all source-$T$-acyclic 1-cells $M \in \sD(C,A)$.  The full
  subcategory of source-$T$-local 1-cells in $\sD(C,A)$ is denoted
  $\sD(C,A)_{\lng \odot T \rng}$.

  A 1-cell $\cell{M'}{B}{C}$ is \emph{target-$T$-acyclic} if $T \odot M
  = 0$.  A 1-cell $\cell{N'}{B}{C}$ is \emph{target-$T$-local} if
  \[
  \sD(B,C)[M',N']_* = 0
  \]
  for all target-$T$-acyclic 1-cells $M' \in
  \sD(B,C)$.  The full subcategory of target-$T$-local 1-cells in
  $\sD(B,C)$ is denoted $\sD(B,C)_{\lng T\odot \rng}$.
\end{defn}

Baker and Lazarev describe the following in the context of
spectra, but their methods generalize to our setting.  The key
observation is that for any 1-cell $P$ whose source is $A$,
$\shom{T}{P}$ is target-$T$-local.  Likewise, if $P'$ is any
1-cell whose target is $B$, $\thom{T}{P'}$ is source-$T$-local.

\begin{prop}[\cite{baker2004thc}]
  \label{prop:general-factorization}
  Let $\cell{T}{A}{B}$ be a 1-cell in $\sD(A,B)$.  The adjunctions induced
  by $T$ factor through the $T$-local pseudofunctors, which is to say that
  we have the following diagrams of adjoint transformations:
  \[\cmxymat{
    {\yda[\sD]{B}} \ar@<.4ex>[rr]^-{T \odot -} \ar@<.4ex>[dr] 
    & 
    & \ar@<.4ex>[ll]^-{\shom{T}{-}} \ar@<.4ex>[dl]^-{\shom{T}{-}} {\yda[\sD]{A}}\\
    & {\yda[\sD]{B}}_{\lng T \odot \rng} \ar@<.4ex>[ul] \ar@<.4ex>[ur]^-{T \odot -} 
    & \\
  }
  \hspace{2pc}
  \cmxymat{
    {\oyda[\sD]{A}} \ar@<.4ex>[rr]^-{- \odot T} \ar@<.4ex>[dr] 
    & 
    & \ar@<.4ex>[ll]^-{\thom{T}{-}} \ar@<.4ex>[dl]^-{\thom{T}{-}} {\oyda[\sD]{B}}\\
    & {\oyda[\sD]{A}}_{\lng \odot T \rng} \ar@<.4ex>[ul]
    \ar@<.4ex>[ur]^-{- \odot T} & \\
  }\]
\end{prop}

\begin{prop}[\cite{baker2004thc}]\label{prop:BLfactorization}
  If a 1-cell $T \in \sD(A,B)$ is right-dualizable and the unit
  map induces an isomorphism $A \iso \thomshort{B}{T}{A}{T}{A},$
  then the induced adjoint pair is an equivalence
  \[
  \yda[\sD]{B}_{\lng T\odot \rng} \hty \yda[\sD]{A}.
  \] 
  We have a corresponding statement for the case of
  left-dualizability.
\end{prop}
\begin{proof}
  It follows immediately from \cref{lem:dlty-eval} that the evaluation
  map 
  \[
  T \odot \shom{T}{M} \to M
  \]
  is an isomorphism for all $M$; this is the counit of the adjunction.
  Moreover, $\shom{T}{-}$ takes values in the target-$T$-local category and
  hence the fact that the unit of the adjunction is an isomorphism follows
  from the fact that the counit is so.
\end{proof}

\begin{cor}\label{cor:faithful-invertible}
  Let $\cell{T}{A}{B}$ be as in \cref{prop:BLfactorization}.  Then $T$
  is target-faithful (\cref{defn:odD}) if and only if localization
  induces an equivalence between the category of 1-cells
  $\scr{B}(B,B)$ and the target-$T$-local subcategory.  In this case
  each of the three adjoint pairs of \fref{prop:general-factorization}
  (at left) is an equivalence.  We have a corresponding statement for
  the source-faithful case.
\end{cor}
\begin{proof}
  If $T$ is target-faithful, then all target-$T$-acyclics are trivial, and
  therefore 
  \[
  \yda[\sD]{B} \hty \yda[\sD]{B}_{\lng T \odot \rng}.
  \]
  The statement then follows from \cref{prop:BLfactorization}.  Conversely,
  if localization induces an equivalence
  $\sD(B,B) \hty \sD(B,B)_{\lng T \odot \rng}$ then the proof of
  \cref{prop:BLfactorization} shows that $T$ is invertible and hence
  target-faithful.
\end{proof}

\begin{defn}\label{defn:faithful-local}
  We say that $B$ is \emph{strongly target-$T$-local} if localization
  induces an equivalence
  \[
  \sD(B,B) \hty \sD(B,B)_{\lng T \odot \rng}.
  \]  
  Likewise, we say that $A$ is \emph{strongly source-$T$-local} if
  localization induces an equivalence
  \[
  \sD(A,A) \hty \sD(A,A)_{\lng \odot T \rng}.
  \]
\end{defn}

Combining the previous results yields a characterization of invertible
objects in triangulated Eilenberg-Watts bicategories.  Applying this
to the case $T = \fl A$ immediately gives the characterization of
Azumaya objects in \cref{thm:Azumaya-char-triang}.
\begin{samepage}
\begin{prop}\label{thm:InvertibilityChar}
  Let $\cell{T}{A}{B}$ be a 1-cell in $\sD$.  The following are equivalent:
  \begin{enumerate}
  \item $T$ is invertible.
  \item
    \begin{enumerate}
    \item $T$ is right-dualizable.
    \item The unit induces $A \iso \thom{T}{T}$.
    \item $B$ is strongly target-$T$-local.
    \end{enumerate}
  \item
    \begin{enumerate}
    \item $T$ is left-dualizable.
    \item The unit induces $B \iso \shom{T}{T}$.
    \item $A$ is strongly source-$T$-local.
    \end{enumerate}
  \end{enumerate}
\end{prop}
\end{samepage}

\begin{rmk}\label{rmk:faithfulness-cond}
  We have the following chain of implications for $\cell{T}{A}{B}$,
  and a similar chain for source-faithful and source-local conditions:
 \[
  T \text{ is target-faithful. } \Rightarrow B \text{ is strongly
    target-$T$-local. } \Rightarrow B \text{ is target-$T$-local.}
  \]
  \Cref{cor:faithful-invertible} shows that the first implication is
  an equivalence under the conditions of \cref{prop:BLfactorization},
  and therefore condition (\emph{c}) in each of the statements of
  \cref{thm:InvertibilityChar} can be replaced by the corresponding
  faithfulness condition on $T$.
\end{rmk}

\subsection{Application to tilting theory}
The work in this section allows us to give a unified proof of results
from the tilting theory of \cite{rickard1989mtd} and \cite{SS03Stable}:
\begin{proof}[Proof of \ref{prop:RiCor} and \ref{prop:RiSpectra}]
  Let $\wt{T}$ denote $T$ regarded as a bimodule over $A = \thom{T}{T}$.
  Since $T$ is right-dualizable, $\wt{T}$ is right-dualizable in
  $\sD_k(A,B).$ Moreover, $A \hty \thom{\wt{T}}{\wt{T}}.$
  
  By \fref{rmk:k-alg-trbicat,prop:homotopy-bicat-applications},
  $\sD_k$ is a triangulated Eilenberg-Watts bicategory.  Since $k$ is the unit of
  $\sD_k$, the 1-cells of $\sD_k(k,B)$ are jointly target-faithful
  (\cref{defn:odD}).  Since $T$ generates $\sD_k(k,B),$
  \cref{lem:odDlem} shows that $T$ is target-faithful and this means
  that $\wt{T}$ is also target-faithful.  The result then follows from
  \cref{cor:faithful-invertible,thm:InvertibilityChar}.
\end{proof}

\section{Homotopical Brauer Groups}
\label{sec:homotopical-brauer-groups}

This section describes homotopical Brauer groups for rings and ring
spectra; these constitute our main applications of the preceding
theory.  We also give explicit comparisons between the Brauer group as
characterized by \fref{thm:AzumayaEquiv,thm:Azumaya-char-triang} and
as it appears in related work on Brauer groups in homotopical
settings.  We begin with the derived Brauer group of a ring and then
address the Brauer groups of ring spectra.

\subsection{The derived Brauer group of a ring}

\begin{defn}
  We refer to $Br(\sD_k)$ as the \emph{derived Brauer group} of $k$, to
  distinguish it from the generally different classical Brauer group, $Br(k)$.
\end{defn}

To\"en \cite{To12Derived} introduces the notion of \emph{derived
  Azumaya algebras} (the two conditions appearing in
\fref{prop:comparison-toen}), and describes a Brauer group formed by
Eilenberg-Watts equivalence classes of such.  We show that the notion
of derived Azumaya algebra is equivalent to the notion of Azumaya
object in the derived category, and therefore the resulting Brauer
groups are isomorphic.
\begin{prop}\label{prop:comparison-toen}
  Let $k$ be a graded commutative ring and let $A$ be a $k$-algebra,
  regarded as a left module over $k$ and a right module over the
  enveloping algebra $A^e$.  Then $A$ is an Azumaya object of $\sD_k$
  if and only if the following two conditions hold:
  \begin{enumerate}
  \item The underlying $k$-module of $A$ is a compact generator of the
    triangulated category $\sD_k(k,k)$.
  \item The map $\mu\cn A^e \to F_k(A,A)$ is an equivalence.
  \end{enumerate}
\end{prop}
\begin{proof}
  To\"en \cite[2.8]{To12Derived} shows that the conditions above are
  equivalent to the condition that $A^e$ be Eilenberg-Watts equivalent to $k$
  in the bicategory $\sD_k$, and this is one of the conditions
  appearing in the characterization \fref{thm:AzumayaEquiv}.
\end{proof}

Combining this with \cite[2.12]{To12Derived} we have the following:
\begin{prop}\label{prop:derived-brauer-group-field}
  Let $k$ be a field.  Then $Br(\sD_k)$ is isomorphic to the classical Brauer group $Br(k)$.
\end{prop}

\noindent
For $k$ not a field, this result does not generalize,
although there is of course always a homomorphism
\[
Br(k) \to Br(\sD_k).
\]

\subsection{The Brauer group of a ring spectrum}

When $k$ is a commutative ring spectrum Baker, Richter and Szymik have
introduced and studied the notion of topological Azumaya $k$-algebra
\cite{BRS10Brauer}.  Their definition (the three conditions of
\fref{prop:comparison-baker-richter} below) makes sense in any modern
monoidal model category of spectra, and we show it is equivalent to
the definition of Azumaya used here.  Since the various modern
categories of spectra are all strong monoidal Quillen equivalent
\cite{SS03Equivalences}, the results of homotopical
Brauer theory transfer between any of them. 

Recall that $\cell{\fl A}{k}{A^e}$ denotes $A$ regarded as a
right module over $A^e$, and ${_k}A_k = \fl{A}_k$ is the underlying
$(k,k)$-module of $A$.
\begin{prop}\label{prop:comparison-baker-richter}
  Let $k$ be a commutative ring spectrum such that $\sD_k$ has the
  structure of a triangulated Eilenberg-Watts bicategory.  Let $A$ be
  a $k$-algebra, regarded as a left module over $k$ and a right module
  over the enveloping algebra $A^e$.  Then $A$ is an Azumaya object in
  $\sD_k$ if and only if the following conditions hold:
  \begin{enumerate}
  \item $A$ is dualizable in the homotopy category of $k$-modules.
  \item The map $A^e \to F_k(A,A)$ is a weak equivalence of bimodules
    over $A^e$.
  \item $A$ is faithful as a $k$-module.
  \end{enumerate}
\end{prop}
\begin{proof}
  Baker, Richter, and Szymik show that these conditions imply that
  $\fl A$ is central and separable over $k$ \cite[1.3,
  1.4]{BRS10Brauer}.  They are working with $S$-algebras, but those
  results apply equally well in any modern category of spectra.  Now
  if $\fl{A}_k = \fl A \odot A_k$ is faithful, then $\fl A$ must be
  target-faithful (\fref{defn:odD}).  \fref{rmk:faithfulness-cond}
  explains that under these circumstances target-faithfulness of $\fl
  A$ is equivalent to $A^e$ being strongly target-$\fl A$-local, and
  thus $A$ is Azumaya by \fref{thm:InvertibilityChar}.

  For the converse, we again refer to \fref{thm:InvertibilityChar}.
  \Fref{lem:duality-composition} shows that $\fl A$ being
  left-dualizable implies $\fl{A}_k$ is left-dualizable, and under
  these circumstances $k$ being strongly source-$\fl A$-local is
  equivalent to $\fl{A}$ being source-faithful
  (\fref{rmk:faithfulness-cond}).  Since restriction of scalars along
  the unit map $k \to A^e$ is faithful, the composite $\fl{A}_k$ must
  also be faithful.
\end{proof}

One minor subtlety remains in the choice of equivalence relation with
which we form a Brauer group.  There are two standard choices,
Eilenberg-Watts equivalence and Brauer equivalence.  We finish by showing
that the two resulting groups are isomorphic.  This is a classical
fact for discrete rings which holds for formal reasons.
\begin{defn}[Brauer equivalence]\label{defn:brauer-equiv}
  Two $k$-algebras $A_1$ and $A_2$ are called (topological) Brauer equivalent if
  there are faithful, dualizable, cofibrant $k$-modules $M_1$ and
  $M_2$ such that
  \[
  A_1 \sma_k F_k(M_1, M_1) \hty A_2 \sma F_k(M_2,M_2)
  \]
  as $k$-algebras.
\end{defn}
\begin{prop}[{\cite[2.4]{BRS10Brauer}}]
  The collection of Azumaya $k$-algebras modulo Brauer
  equivalence forms a group.
\end{prop}
\begin{lem}
  The group of Brauer equivalence classes of Azumaya $k$-algebras is
  isomorphic to the group of Eilenberg-Watts equivalence classes of
  such algebras.
\end{lem}
\begin{proof}
  For a dualizable $k$-module $M$, $F_k(M,M)$ is always
  Eilenberg-Watts equivalent to $k$.  So if $A_1$ is Brauer equivalent
  to $A_2$ then it is Eilenberg-Watts equivalent to $A_2$.  Therefore Brauer
  equivalence implies Eilenberg-Watts equivalence and there is a
  surjective group homomorphism from Azumaya algebras modulo Brauer
  equivalence to Azumaya algebras modulo Eilenberg-Watts equivalence.
  But if $A$ is Eilenberg-Watts equivalent to the ground ring, $k$,
  then there is an invertible bimodule $M$ giving the equivalence, and
  $A \hty F_k(M,M)$, so $A$ is also Brauer equivalent to $k$.  Thus
  this homomorphism is an isomorphism.
\end{proof}

\bibliographystyle{amsalpha}
\bibliography{final-version.bbl}

\end{document}